\documentclass{amsart}
\usepackage{amsmath,amssymb,amscd,graphicx,psfrag,epsf,amsthm}
\usepackage[all]{xy}

\newcommand{\F}{{\mathbb {F}}}

\newcommand{\cC}{{\mathcal {C}}}

\newcommand{\cE}{{\mathcal {E}}}

\newcommand{\cD}{{\mathcal {D}}}

\newtheorem{thm}{Theorem}
\newtheorem{lemma}[thm]{Lemma}

\newtheorem{prop}[thm]{Proposition}

\newtheorem{conj}[thm]{Conjecture}
\newtheorem*{thm*}{Theorem}
\newtheorem*{mainthm}{Main Theorem}

\theoremstyle{definition}
\newtheorem{defn}[thm]{Definition}

\newtheorem{say}[thm]{}

\newtheorem{notation}[thm]{Notation}   
  
\newtheorem{defn-thm}[thm]{Definition-Theorem}  
  
\newtheorem{rem}[thm]{Remark} 
\newtheorem{remarks}[thm]{Remarks} 
 
\newtheorem*{notation*}{Notation}

\theoremstyle{remark}

\renewcommand{\o}[0]{{\mathcal O}}  
\newcommand{\z}[0]{{\mathbb Z}}

\renewcommand{\r}[0]{{\mathbb R}}

\newcommand{\p}[0]{{\mathbb P}}

\newcommand{\q}[0]{{\mathbb Q}}
\newcommand{\map}[0]{\dasharrow}

\newcommand{\rat}[0]{\operatorname{RatCurves}^n}
\newcommand{\pic}[0]{\operatorname{Pic}}

\newcommand{\codim}[0]{\operatorname{codim}}

\newcommand{\locus}[0]{\operatorname{locus}}

\newcommand{\chow}[0]{\operatorname{Chow}}

\def\into{\DOTSB\lhook\joinrel\rightarrow}

\numberwithin{equation}{section}

\begin{document}

\title[Identifying quadric bundles]{Identifying quadric bundle structures on complex projective varieties}

\author{Carolina Araujo}

\address{Carolina Araujo: \sf IMPA, Estrada Dona Castorina 110, Rio de
  Janeiro, 22460-320, Brazil} 
\email{caraujo@impa.br}


\maketitle

\begin{abstract}
In this paper we characterize smooth complex projective varieties that admit a quadric bundle structure 
on some dense open subset in terms of the geometry of  certain families of rational curves. 
\end{abstract}


\section{Introduction}\label{introduction}

Let $X$ be a smooth complex projective variety, and assume that $X$
is uniruled, i.e., there exists a rational curve through every point of $X$.
Then there exists an irreducible family of rational curves that sweeps out a dense open subset
of $X$. Such a family is called a \emph{covering} family of rational curves on $X$.  
A covering family of rational curves reflects a lot of 
the geometry of the ambient variety $X$, specially when this family is ``minimal'' in some sense. 
There are several different notions of ``minimality''.
In this paper we consider 
\emph{unsplit} and \emph{minimal} covering families. 
A covering family $H$ of rational curves on $X$ is called \emph{minimal} if, for a general point $x\in X$,
the subfamily $H_x$ of $H$ parametrizing curves through $x$ is proper. 
It is called \emph{unsplit} if $H$ itself is proper.
Unsplit covering families of rational curves are very powerful tools in studying the geometry of  uniruled varieties.
However, not every uniruled variety carries such a family.
Think, for instance, of a smooth quartic threefold. The family of conics on it is minimal but not unsplit.
On the other hand, every uniruled complex projective  variety $X$ carries a minimal  covering family of rational curves.
For instance, one may fix an ample line bundle on $X$ and take $H$ to be a 
covering family of rational curves having minimal degree with respect to this fixed line bundle.

Fix $H$ a minimal covering family of rational curves on $X$.
Given a general point $x\in X$, let $\cC_x\subset \p(T_xX)$ denote the closed subset of the 
projectivized tangent space at $x$ consisting of 
tangent directions at $x$ to rational curves parametrized by $H_x$ (see Definition~\ref{def_Cx}).
This is called the \emph{variety of minimal rational tangents} at $x$ associated to the family $H$.
It was first introduced in \cite{mori79}, where Mori used it to prove the Hartshorne conjecture, and
it has found many applications within the theory of uniruled varieties, specially Fano varieties
(see \cite{hwang} for a survey).
The variety of minimal rational tangents comes with a natural embedding into $\p(T_xX)$.
The general philosophy is that many properties of $X$  can be detected by studying this projective embedding.
For instance, Theorem~\ref{thm_CMSB} below establishes that $\cC_x= \p(T_xX)$ if and only if $X\cong \p^n$.

\begin{thm}[{\cite{CMSB} - see also \cite{kebekus_on_CMSB}}] \label{thm_CMSB}
  Let $X$ be a smooth complex projective $n$-dimensional variety, $H$ a minimal
  covering family of rational curves on $X$, and $\cC_x\subset
  \p(T_xX)$ the associated variety of minimal rational tangents
  at $x\in X$. Suppose that  $\cC_x=\p(T_xX)$ for a general point $x\in X$.
  Then $X\cong \p^n$ and under this isomorphism $H$ corresponds to the family of lines on $\p^n$.
\end{thm}

If $X\not\cong \p^n$, then it is not possible to completely recover $X$
from the embedding $\cC_x\into \p(T_xX)$. 
Indeed, if $X'$ is the blow up of $X$ at a point, then $H$ induces a
minimal covering family of rational curves on $X'$ whose 
variety of minimal rational tangents at a general point is 
projectively isomorphic to $\cC_x$.
On the other hand, if one imposes some extra condition on $X$, for instance the 
Picard number being $1$, then in some cases $X$ may be fully recovered 
from its variety of minimal rational tangents at a general point. 
This is illustrated by the next theorem.

\begin{thm}[{\cite[Main Theorem]{mok_05}}] \label{thm_mok}
  Let $X$ and $Y$ be smooth complex projective Fano varieties with Picard number $1$, and suppose that 
  $Y$ is either a Hermitian symmetric space or a contact homogeneous manifold. 
  Let $H_X$ and $H_Y$ be minimal covering families of rational curves on $X$ and $Y$ respectively, 
  and $\cC_x\subset \p(T_xX)$ and $\cD_y\subset \p(T_yY)$ the associated varieties of minimal rational tangents
  at general points $x\in X$ and $y\in Y$, respectively. 
  Suppose that there is an isomorphism $\p(T_xX)\cong \p(T_yY)$ inducing an isomorphism 
  $\cC_x\cong \cD_y$. 
  Then  $X\cong Y$, and under this isomorphism curves from $H_X$ correspond to curves from $H_Y$.
\end{thm}

Notice that this theorem applies in particular to a smooth hyperquadric $Y=Q_n\subset \p^{n+1}$ with $n\geq 3$.
In this case the only minimal covering family of rational curves on $Y$ is the family of lines, and 
the associated variety of minimal rational tangents at any $y\in Y$ is a 
positive dimensional hyperquadric in $\p(T_yY)$.

When the Picard number of $X$ is not $1$, one can still use $\cC_x$ to get information about $X$.
The following is a structure theorem for varieties $X$ for which  $\cC_x$  is a linear subspace of $\p(T_xX)$
for general $x\in X$.

\begin{thm}[{\cite[Theorems 1.1 and 3.4]{artigo_tese}}]\label{pn_bundle_in_codim_1} 
  Let $X$ be a smooth complex projective variety, $H$ a minimal
  covering family of rational curves on $X$, and $\cC_x\subset
  \p(T_xX)$ the associated variety of minimal rational tangents
  at $x\in X$. Suppose that for a general point $x\in X$, $\cC_x$ is a
  $d$-dimensional linear subspace of $\p(T_xX)$.
  
  Then there exists a dense open subset $X^\circ\subset X$, a smooth quasi-projective variety $Y^\circ$ and a
  $\p^{d+1}$-bundle $\pi^\circ:X^\circ\to Y^\circ$ such that
  every  curve from  $H$ meeting $X^\circ$ corresponds to a line on a fiber of $\pi^\circ$. 
  If $H$ is unsplit, then one may take $X^\circ$ such that
  $\codim_X(X\setminus X^\circ)\geq 2$.
\end{thm}

In this paper we go one step further. We provide the following 
structure theorem for varieties $X$ carrying an unsplit covering family of rational curves for which
$\cC_x$ is a positive dimensional hyperquadric in its linear span in $\p(T_xX)$ for general $x\in X$. 
By a \emph{quadric bundle} we mean a flat projective morphism between quasi-projective varieties whose fibers are all isomorphic to 
 irreducible and reduced (but not necessarily smooth) hyperquadrics.

\begin{mainthm}
  Let $X$ be a smooth complex projective variety admitting 
  an unsplit  covering family of rational curves $H$.
  Let  $\cC_x\subset \p(T_xX)$ denote the variety of minimal rational tangents  
  associated to $H$ at a general point $x\in X$. 
  Suppose $\cC_x$ is a positive dimensional hyperquadric in its linear span in $\p(T_xX)$.

  Then there exists an open subset $X^\circ\subset X$, with $\codim_X(X\setminus X^\circ)\geq 2$,
  a smooth quasi-projective variety $Y^\circ$ and a quadric bundle $\pi^\circ:X^\circ\to Y^\circ$ such that
  every curve from $H$
  meeting $X^\circ$ corresponds to a line on a fiber of $\pi^\circ$. 
\end{mainthm}

\begin{remarks} {\bf (1)}
The unsplitness assumption in the Main Theorem cannot be dropped. 
Indeed, start with a smooth complex projective variety $Z$ admitting a quadric bundle structure $\pi:Z\to Y$, 
and let $\sigma\subset Z$ be a multi-section of $\pi$. 
Then let $X$ be the blowup of $Z$ along $\sigma$.

\noindent {\bf (2)} 
In general one should not expect the quadric bundle $\pi^\circ:X^\circ\to Y^\circ$ obtained in the Main Theorem
to extend to a quadric bundle in the whole of $X$. For instance, given $n\geq 6$, let $X$ be the complete intersection in $\p^n\times \p^n$ of two general divisors 
of type $(1,1)$, and one general divisor of type $(0,2)$. Let $\pi:X\to \p^n$ be the first projection. The general fiber of $\pi$
is isomorphic to a ($n-3$)-dimensional hyperquadric, but there are fibers of higher dimension.

\noindent {\bf (3)} 
In the special case when $H$ is a family of lines under some projective embedding of $X$, and $\dim \cC_x > \left[\frac{\dim X}{2}\right]$, 
this result follows from \cite[Theorem 3.1]{beltrametti_ionescu_q-bdles}.
\end{remarks}




This paper is organized as follows. 
In Section~\ref{section_minimal_rat_curves} we gather results
about minimal and unsplit covering families of rational curves and their
rationally connected quotients.  In Section~\ref{section_Cx}
we discuss the variety of minimal rational tangents associated 
to a minimal covering family of rational curves and the distribution defined by its linear span.
In Section~\ref{section_proof}, we prove the Main Theorem. 

\begin{notation*}
  Throughout this paper we work over the field of complex
  numbers.  
  By a rational curve we always mean a projective rational curve.
  If $V$ is a complex vector space,
  we denote by $\p(V)$ the projective space of $1$-dimensional linear subspaces of $V$.
  Let $X$ be a quasi-projective variety. 
  By a general point of $X$, we mean a point in some dense open subset of $X$.
  If $f:X\to T$ is a proper morphism onto another quasi-projective variety and $t\in T$, we denote the fiber of $f$ over $t$ by $X_t$.
  We denote by  $\rho(X/T)$ the relative Picard number of $X$ over $T$. 
  This is the dimension of the vector space $N_1(X/T)$ of $1$-cycles on $X$ with real coefficients 
  generated by irreducible curves contracted by $f$, modulo
  numerical equivalence.
\end{notation*}


\section{Families of rational curves and their
rationally connected quotients}\label{section_minimal_rat_curves}

Let $X$ be a complex projective variety.
There is a scheme $\rat(X)$ parametrizing rational curves on $X$.
This scheme is constructed as the normalization of a certain 
subscheme of the Chow scheme $\chow(X)$ parametrizing effective $1$-cycles on $X$.
More generally, given a proper morphism $f:X\to Y$ between complex quasi-projective varieties,
there is a $Y$-scheme $\rat(X/Y)$ parametrizing rational curves on $X$ contained on fibers of $f$.
As before, $\rat(X/Y)$ is  constructed as the normalization of a certain 
subscheme of the Chow scheme $\chow(X/Y)$ parametrizing effective $1$-cycles on $X$ contracted by $f$.
We refer to \cite[Chapters I and II]{kollar} for details, constructions and proofs.
See also \cite{debarre}.

\begin{defn}\label{def_H}
Let $f:X\to Y$ be a proper morphism  between complex quasi-projective varieties.
By a \emph{family of rational curves on $X$ over $Y$} we mean  an
irreducible component $H$ of $\rat(X/Y)$. 
In particular $H$ is a normal quasi-projective variety.
The universal properties of $\rat(X/Y)$ yield  universal family morphisms 
\[
\begin{CD}
  U @>\eta >> X, \\
  @V{p} VV \\
  H
\end{CD}
\]
where $U$ is a normal quasi-projective variety, $p:U\to H$ is a $\p^1$-bundle, and $\eta:U\to X$
is the evaluation morphism. 
We say that $H$ is a \emph{covering family} if the image of $\eta$ is dense in $X$.

Now assume that $Y$ is a point and $X$ is a complex projective variety, and hence $\rat(X/Y)=\rat(X)$.
Let $H$ be  a  covering family of rational curves on $X$.
Given a point $x\in X$, we denote by $H_x$  the 
subscheme of $H$ parametrizing rational curves passing through $x$, and set
$\locus(H_x):=\eta\big(p^{-1}(H_x)\big)$ (with the reduced scheme structure).
We say that $H$ is \emph{unsplit} if it is proper. 
We say that it is \emph{minimal} if, for a 
general point $x\in X$, $H_x$ is proper.
\end{defn}

\begin{rem}\label{restriction_of_H}
Let $f:X\to Y$ be a proper morphism  between complex quasi-projective varieties, and $y\in Y$ a point.
Then there exists a natural map 
$$
\iota: \rat(X_y)\to \rat(X/Y).
$$
Let $H'$ be a family of rational curves on $X_y$, and let $H$ be the irreducible component of $\rat(X/Y)$ 
containing $\iota\big(H'\big)$.
Then it may happen that $\iota^{-1}(H)\neq H'$. 
For instance, let $f:X\to Y$ be a family of quadric surfaces in $\p^3$,
$X_y$ a smooth fiber, and $H'$ one of the two families of lines on $X_y$.  If 
we choose $f:X\to Y$ suitably, then $\iota^{-1}(H)$ is the union of the two families of lines on $X_y$.

Let $C, C' \subset X_y$ be two curves parametrized by points of $H$. Then they are numerically equivalent on $X$.
If the restriction morphism $\pic(X)\otimes_\z \q \to \pic(X_y)\otimes_\z \q$ is surjective, then  they are also numerically equivalent on $X_y$.

This observation allows us to conclude in certain cases that  $\iota^{-1}(H)= H'$.
For instance, let   $f:X\to Y$ be a proper morphism  between complex quasi-projective varieties, 
and suppose that $X_y$ is isomorphic to a smooth hyperquadric 
$Q_n\subset \p^{n+1}$, with $n\geq 3$.
Let $H'$ be the family of rational curves on $X_y$ corresponding 
to lines on $Q_n$ under the isomorphism $X_y\cong Q_n$. 
Then $\pic(X_y)\cong \z$, and thus
the morphism $\pic(X)\otimes_\z \q\to \pic(X_y)\otimes_\z \q$ is surjective. 
Any curve on $Q_n$ that is numerically equivalent to a line is itself a line.
Moreover,  there is a unique family of lines on $Q_n$.
Hence $\iota^{-1}(H)= H'$.
\end{rem}

\begin{defn}[{{\bf The $H$-rationally connected quotient}}] \label{defn_HRQ}
Let $X$ be a normal complex projective variety.
Suppose that $X$ admits an unsplit  covering family of rational curves $H$.
We define an equivalence relation $\sim_H$ on $X$ as follows: 
$$
x\sim_H y \ \iff \ x \text{ and } y \text{ can be connected by a chain of curves from } H.
$$
By  \cite[IV.4.16]{kollar} (see also \cite{campana}), there exists a
proper surjective morphism $\pi^\circ:X^\circ \to Y^\circ$ from a
dense open subset of $X$ onto a normal quasi-projective variety whose fibers are 
$\sim_H$-equivalence classes.  We call this map the 
\emph{$H$-rationally connected quotient of $X$}.  
When $Y^\circ$
is a point we say that $X$ is $H$-rationally connected.
\end{defn}

\begin{remarks} \label{rems_HRQ} 
{\bf (1)}  If $X$ is smooth, then, by \cite[Lemma 2.2]{ADK_Beauville}, there is a morphism $\pi^\circ:X^\circ \to Y^\circ$
as above with the additional properties that $\codim_X(X\setminus X^\circ)\geq 2$,  $Y^\circ$ is smooth, and 
$\pi^\circ$ is an equidimensional proper morphism with irreducible and reduced fibers
(see also \cite[Proposition 1]{BCD_extremal_rays} and \cite[3.1, 3.2]{andreatta_wisniewski}).

\smallskip

\noindent {\bf (2)}  If $X$ is smooth and $H_x$ is irreducible for general $x\in X$, then 
the general fiber of  $\pi^\circ$ is a smooth Fano variety with Picard number $1$
by \cite[Propopsition 2.3]{ADK_Beauville}.
\end{remarks}

We end this section with a simple but useful observation.

\begin{prop}\label{rho(X/Y)=1}
Let $f:X\to C$ be a proper morphism from a normal quasi-projective variety onto a smooth quasi-projective curve.
Suppose that the general fiber of $f$ is a smooth Fano variety with Picard number $1$, and that every fiber
of $f$ is irreducibe.
Then $\rho(X/C)=1$.
\end{prop}

\begin{proof}
Let $H$ be  a  covering family of rational curves on $X$ over $C$.
Let $C^\circ\subset C$ be a dense open subset such that, for every $c\in C^\circ$,
$X_c$ is a smooth Fano variety with Picard number $1$, and
$X_c$ contains a curve parametrized by some point of $H$.
We note that, for every $c\in C^\circ$, $\pic(X_c)\cong \z$ and $N_1(X_c)=\r\cdot [\ell]$, where $[\ell]\in N_1(X_c)$
is the class of a curve  $\ell \subset X_c$ parametrized by a point of $H$.
Set $X^\circ =f^{-1}(C^\circ)$, and $f^\circ=f|_{X^\circ}: X^\circ \to C^\circ$.

 First we observe that $\pic(X^\circ)=(f^\circ)^*\pic(C^\circ) \oplus \z\cdot L^\circ$, where 
 $L^\circ\in \pic(X^\circ)$ is $(f^\circ)$-ample.
 Indeed, let $F$ denote a general fiber of $f^\circ$, 
 $r:\pic(X^\circ)\to \pic(F)\cong \z$ the natural restriction map,  
 and suppose that $M\in \pic(X^\circ)$ is such that $r(M)=0$.
 Then $M\cdot \ell=0$ for every curve $\ell\subset X^\circ$ parametrized by some point of $H$.
 Thus $M|_{X_c}\cong \o_{X_c}$ for every $c\in C^\circ$.
 Therefore $(f^\circ)_*M$ is an invertible sheaf on $C^\circ$ (see for instance \cite[III.12.9]{hartshorne}), 
 and there is an isomorphism $(f^\circ)^*(f^\circ)_*M\to M$, i.e., $M\in (f^\circ)^*\pic(C^\circ)$.
 
Let $L$ be an element of $\pic(X)$ extending $L^\circ$.  Since all fibers of $f$ are irreducible,
$\pic(X)\otimes_\z \q=(f^*\pic(C)\otimes_\z \q) \oplus \q\cdot L$ by \cite[II.6.5]{hartshorne}.
From this it follows that $\rho(X/C)=1$.
\end{proof}


\section{The variety of minimal rational tangents}\label{section_Cx}

Let $X$ be a smooth complex projective uniruled variety,  $H$ a minimal covering family of rational curves on $X$, 
and $x\in X$ a general point.

\begin{defn} \label{def_Cx}
We define the map $\tau_x:  H_x  \map \p(T_xX)$
by sending a curve that is smooth at $x$ to its tangent direction at $x$.
We define the \emph{variety of minimal rational tangents} associated to $H$ at $x$ to be the closure of the image of $\tau_x$
in $\p(T_xX)$, and denote it by $\cC_x$.
Notice that $\cC_x$ may be reducible.
\end{defn}

\begin{remarks} \label{rems_Cx} 
{\bf (1)}  Let $\tilde H_x$ be the normalization of $H_x$.
By \cite[II.1.7, II.2.16]{kollar}, $\tilde H_x$ is a smooth projective variety.
Moreover, the map $\tilde \tau_x: \tilde H_x \to \cC_x$ induced by $\tau_x:  H_x  \map \cC_x$ is the normalization
morphism by \cite{kebekus} and \cite{hwang_mok_birationality}. 
This implies in particular that $\cC_x$ cannot be a singular cone in $\p(T_xX)$
(see \cite[Lemma 4.3]{artigo_tese}).

\smallskip

\noindent {\bf (2)}  If $\cC_x$ is a union of linear subspaces of $\p(T_xX)$, then
the intersection of any two irreducible components of $\cC_x$ is empty by  \cite[Lemma 2.8]{ADK_Beauville}
(see also \cite[Proposition 2.2]{hwang_b2=b4=1}). 
\end{remarks}

In the next section, we investigate varieties $X$ for which $\cC_x$ is 
a positive dimensional hyperquadric in its linear span in $\p(T_xX)$.
Remarks~\ref{rems_Cx}  above imply that in this case $\cC_x$ is a 
smooth hyperquadric in its linear span in $\p(T_xX)$.

Our next goal is to explain how the  variety of minimal rational tangents defines 
a distribution on $X$, and provide conditions under which this distribution is integrable. 
First let us recall the definition of distribution and Frobenius'  criterion for integrability.

\begin{say}[{{\bf Distributions}}]\label{distributions}
Let $X$ be a smooth complex projective variety.
A \emph{distribution} on $X$ is a subbundle of the tangent bundle of a 
dense open subset $X^\circ$ of $X$, $E^\circ\into X^\circ$.
We regard two distributions as being the same if they agree on the open subset
where both are defined. 
The subbundle $E^\circ\into T_{X^\circ}$  can be extended to a subsheaf $E$ of $T_X$. 
By abuse of notation we denote this distribution by $E\into T_X$.

We say that the distribution $E\into T_X$ is \emph{integrable} if through every point $x\in X^\circ$
there is a complex analytic manifold $M\subset X^\circ$  such that $E_y = T_yM$ 
for every $y\in M$. 

The Lie bracket defines a section of $Hom(\wedge^2 E^\circ,T_{X^\circ}/E^\circ)$,
which we call the \emph{Frobenius bracket tensor of $E$},
$$
[,] \ : \  \wedge^2 E^\circ \ \longrightarrow \ T_{X^\circ}/E^\circ.
$$
By Frobenius' Theorem, the distribution $E\into T_X$ is integrable if and only if 
the Frobenius bracket tensor of $E$ vanishes identically.
\end{say}

Now we fix notation and gather assumptions that will be used for the rest of the paper.

\begin{notation} \label{assumptions}
 {\bf (1)} Let $X$ be a smooth complex projective variety, and $H$ a minimal covering family of rational curves on $X$.
Let $\cC_x$ be the  variety of minimal rational tangents associated to $H$ at a general point $x\in X$,
and denote by $\cE_x$ its linear span in $\p(T_xX)$.
Then we obtain a distribution $E\into T_X$ by setting 
$E_x$ to be the linear subspace of $T_xX$ corresponding 
to $\cE_x \subset \p(T_xX)$ for general $x\in X$.
We denote by
$$
[,]_x \ : \  \wedge^2 E_x \ \longrightarrow \ T_xX/E_x,
$$
the Frobenius bracket tensor of $E$ at a general point $x$, 
and by $C_x$ the cone in $E_x$ corresponding to $\cC_x\subset \cE_x$.

\smallskip

\noindent {\bf (2)} 
Let $V$ be a complex vector space, and $C\subset V$ a cone such that $\p(C)$ is an algebraic set in $\p(V)$.
We define $\hat C$ to be the cone in $\wedge^2 V$ consisting of elements of the form $u\wedge v$, where 
$u$ and $v$ generate a linear subspace of $V$ that is tangent to the cone $C$ at a smooth point.
Given $\omega\in \wedge^2V^*$, we say that $C$ is \emph{isotropic with respect to $\omega$} if $\omega(\xi)=0$ 
for every $\xi\in \hat C$. 
\end{notation}

\begin{rem} \label{isotropic}
If $V$ and $C$ are as in Notation~\ref{assumptions}.2, $\omega\in \wedge^2V^*$
is nondegenerate (i.e., a symplectic form), and $C$ is isotropic with respect to $\omega$, 
then $\dim C\leq \frac{\dim V}{2}$.
\end{rem}

\begin{prop}[{\cite[Proposition 9]{hwang_mok_hermitian}}] \label{kernel_of_frobenius}
 Let $X$, $H$, $\cC_x$, $\cE_x$,  $E\into T_X$,  $[,]_x$ and $\hat C_x$ be as in Notation~\ref{assumptions}.
Then the kernel of $[,]_x$ contains $\hat C_x$.
\end{prop}

In \cite[Proposition 9]{hwang_mok_hermitian} $X$ is assumed to be a Fano variety, but this condition is not used in the proof.

The next sufficient condition for integrability of $E\into T_X$
is a straightforward consequence of Frobenius' Theorem, Proposition~\ref{kernel_of_frobenius}
and the following linear algebra lemma.

\begin{prop}\label{distribution_is_integrable}
Let $X$, $H$, $\cC_x$, $\cE_x$ and $E\into T_X$ be as in Notation~\ref{assumptions}.
If $\cC_x$ is an irreducible hypersurface  in $\cE_x$ for general $x\in X$, then $E\into T_X$ is integrable.
\end{prop}

\begin{lemma}
Let $V$, $C$ and $\hat C\subset \wedge^2 V$ be as in Notation~\ref{assumptions}.2.
If $C$ is an irreducible nonlinear cone of codimension $1$ in $V$, then $\hat C$ is
nondegenerate in $\wedge^2 V$.
\end{lemma}

\begin{proof}
Suppose that $\hat C$ is degenerate in $\wedge^2 V$.
Then there is a nonzero element $\omega\in \wedge^2V^*$ such that $C$ is isotropic with respect to $\omega$.
Set 
$$
Q=\ker \omega=\big\{v\in V \ \big| \ \omega(v,u)=0 \text{ for every } u\in V\big\},
$$
and consider the natural 
projection $\pi:V\to V/Q$. The form $\omega$ induces a symplectic form $\omega_\pi$ on $V/Q$.
Moreover the cone $\pi(C)$ is isotropic with respect to $\omega_\pi$, and thus 
$\dim \pi(C)\leq \frac{\dim \ V/Q}{2}$ by Remark~\ref{isotropic}.
Since $C$ is nondegenerate, $\pi(C)\neq 0$, and since $\omega_\pi$ is nonzero,  $\pi(C)\neq V/Q$.
Therefore $\pi(C)$ has codimension $1$ in $V/Q$.
So we conclude that $\dim \ V/Q =2$ and $\dim \pi(C)=1$.
Since $\pi(C)$ is an irreducible cone, it must be a $1$-dimensional linear subspace in $V/Q$.
Thus $C$ is a hyperplane in $V$, which contradicts the nonlinearity of $C$.
\end{proof}

When $X$ is a smooth Fano variety with Picard number $1$, we have the following 
necessary condition for integrability of $E\into T_X$.

\begin{prop}[{\cite[Proposition 2]{hwang_stability_of_TX}}] \label{prop_not_integrable}
Let $X$, $H$ and  $\cC_x$ be as in Notation~\ref{assumptions},
and assume moreover that $X$ is a Fano variety with Picard number $1$.
Let $D\into T_X$ be a distribution on $X$ such that $\cC_x\subset \p(D_x)$ for general $x\in X$.
Then $D\into T_X$ is integrable if and only if $D_x=T_xX$ for general $x\in X$.
\end{prop}

\begin{lemma}\label{distribution=HRCQ}
Let $X$, $H$, $\cC_x$, $\cE_x$ and $E\into T_X$ be as in Notation~\ref{assumptions}.
Let $\pi^\circ:X^\circ \to Y^\circ$ be the $H$-rationally connected quotient of $X$,
and assume that the general fiber of $\pi^\circ$ is a Fano variety with Picard number $1$.
If $E\into T_X$ is  integrable, then it coincides with the distribution defined by $\pi^\circ$.
\end{lemma}

\begin{proof}
Let $F$ a general fiber of  $\pi^\circ$, and let $x\in F$ be a general point.
Then $\locus(H_x)\subset F$. 
The projectivized tangent cone to $\locus(H_x)$ at $x$ contains $\cC_x$.
Hence  $\cC_x\subset \p(T_xF)$, and thus $\cE_x\subset  \p(T_xF)$.
Therefore, since the cokernel of $T_F\into T_X|_F$ is torsion free, 
the inclusion $E|_F\into T_X|_F$ factors through an inclusion  $E|_F\into T_F$.
If $E\into T_X$ is an integrable distribution on $X$, then $E|_F\into T_F$ is an integrable distribution on $F$.

By hypothesis $F$ is a smooth Fano variety with Picard number $1$.
Consider the natural map $\iota: \rat(F)\to \rat(X)$. 
There is an irreducible component $H_F$ of $\iota^{-1}(H)$ that is a minimal
covering  family of rational curves on $F$. Moreover, the variety of minimal rational tangents associated  to $H_F$
at a general point $x\in F$ is contained in $\cE_x\subset \p(T_xF)$.
So we can apply Proposition~\ref{prop_not_integrable},  and conclude that if  $E|_F\into T_F$ is integrable, then
$E_x=T_xF$ for $F$ a general fiber of $\pi^\circ$ and $x$ a general point of $F$. I.e., 
$E\into T_X$ coincides with the distribution defined by $\pi^\circ$.
\end{proof}


\section{Quadric bundles}\label{section_proof}

We start this section by recalling the following characterization of quadric bundles due to Fujita.
We say that a quadric bundle  $\pi:X\to Y$ is a \emph{geometric quadric bundle} if there exists a vector bundle 
$V$ on $Y$ such that $X$ embeds into $\p(V)$ over $Y$ as a divisor of relative degree $2$.

\begin{prop}[{\cite[Corollary 5.5]{fujita75}}] \label{fujita}
Let $X$ and $Y$ be irreducible and reduced complex analytic spaces, and $\pi:X\to Y$ a proper and flat morphism 
whose fibers are all irreducible and reduced. 
Suppose that the general fiber of $\pi$ is isomorphic to a hyperquadric $Q_n\subset \p^{n+1}$, and that
there exists a $\pi$-ample line bundle $L$ on $X$ such that 
$L|_{X_t}\cong \o_{\p^{n+1}}(1)|_{Q_n}$ for general $t\in Y$.
Then $\pi$ is a geometric quadric bundle. 
\end{prop}

\begin{prop} \label{generic_q-bdle}
Let $X$ be a smooth complex quasi-projective variety, $B$ a smooth complex quasi-projective curve, and
$\pi:X\to B$ a proper morphism with irreducible fibers.
Suppose that the general fiber of $\pi$ is isomorphic to a hyperquadric $Q_n\subset \p^{n+1}$ with $n\geq 3$.
Then $\pi$ is a geometric quadric bundle. 
\end{prop}

\begin{proof}
We want to apply Proposition~\ref{fujita} to $\pi:X\to B$.
First we notice that, since $X$ and $B$ are smooth and $\pi$ is proper and equidimensional, 
$\pi$ is flat by \cite[6.1.5]{EGA4}.
Moreover, the morphism $\pi$ admits a section $C\subset X$ by \cite{GHS}. Thus 
its fibers are generically reduced, and hence reduced. 
So, in order to prove the lemma, we only need to produce a $\pi$-ample line bundle $L$ on $X$
whose restriction to a general fiber of $\pi$ is isomorphic to $\o_{\p^{n+1}}(1)|_{Q_n}$.

Let $F$ be a general fiber of $\pi$, and  
$H_F\subset \rat(F)$ the unsplit covering family of rational curves on $F$ corresponding 
to lines on $Q_n$ under the isomorphism $F\cong Q_n$. 
Let $H$ be the irreducible component of $\rat(X/B)$  containing the image of $H_F$
under the natural map $\iota: \rat(F)\to \rat(X/B)$.
Since $\rat(X/B)$ has countably many components and $F$ is a general fiber of $\pi$,
$H$ is a covering family of rational curves on $X$ over $B$.
Moreover, by Remark~\ref{restriction_of_H},  $\iota^{-1}(H)= H_F$.
Consider the universal family morphisms (see Definition~\ref{def_H}):
\[
\begin{CD}
  U @>\eta >> X. \\
  @V{p} VV \\
  H
\end{CD}
\]

Let $D$ be the unique irreducible component of the closure of $\eta\Big(p^{-1}\big(p(\eta^{-1}C)\big)\Big)$ in $X$
that dominates $B$ (with the reduced scheme structure).
By construction, $D$ is a Cartier divisor on $X$, and its restriction to a general  fiber of $\pi$ is a member of $\big|\o_{\p^{n+1}}(1)|_{Q_n}\big|$.
By Proposition~\ref{rho(X/Y)=1}, $\rho(X/B)=1$, and thus $D$ is $\pi$-ample.
\end{proof}

\begin{rem}
The condition $n\geq 3$ in Proposition~\ref{generic_q-bdle} cannot be relaxed.
Indeed,  a smooth quadric surface $Q_2$ can degenerate into any Hirzebruch surface $\F_a$ for $a$ even.
\end{rem}

\begin{proof}[Proof of the Main Theorem]
Let the notation be as in Notation~\ref{assumptions}.
By Remarks~\ref{rems_Cx}, $\cC_x$ is a  positive dimensional smooth hyperquadric in $\cE_x$. 
In particular both $\cC_x$ and $H_x$ are irreducible.

Let $\pi^\circ:X^\circ \to Y^\circ$ be the $H$-rationally connected quotient of $X$, where 
$Y^\circ$ is smooth, $\codim_X(X\setminus X^\circ)\geq 2$,   and 
$\pi^\circ$ is an equidimensional proper morphism with irreducible and reduced fibers
(see Remark~\ref{rems_HRQ}.1).
The general fiber $F$ of $\pi^\circ$ is a smooth Fano variety with Picard number $1$
by Remark~\ref{rems_HRQ}.2. Set $n=\dim F$.

By Proposition~\ref{distribution_is_integrable}, $E\into T_X$ is integrable.
Hence, by Lemma~\ref{distribution=HRCQ}, it coincides with the distribution defined by $\pi^\circ$.
Therefore the general fiber $F$ of $\pi^\circ$ 
admits an unsplit family of rational curves whose associated variety of minimal rational tangents
at a general point $x\in F$ is a positive dimensional smooth hyperquadric in $\p(T_xF)$.
By Theorem~\ref{thm_mok}, $F\cong Q_n$, with $n\geq 3$.

Let $B\subset Y^\circ$ be a curve obtained as a complete intersection of general very ample divisors on $Y^\circ$.
Set $X_B =\pi^{-1}(B)$, and $\pi_B=\pi|_{X_B}: X_B \to B$.
By Bertini's Theorem, both $B$ and $X_B$ are smooth.
Moreover, the general fiber of $\pi_B$ is isomorphic to  $Q_n$, with $n\geq 3$,
and every fiber is irreducible.
Thus $\pi_B$ is a (geometric) quadric bundle by Proposition~\ref{generic_q-bdle}.
Therefore the locus $S$ of $Y^\circ$ over which $\pi^\circ$ is not a quadric bundle 
has codimension  at least $2$. By replacing $Y^\circ$ with $Y^\circ\setminus S$ 
we get the desired statement.
 \end{proof}

\medskip
\noindent {\bf Acknowledgements.}
I would like to thank Giuseppe Borrelli and St\'ephane Druel for many fruitful discussions
about some of the topics that appear in this paper.
Partial financial support has been provided by CNPq -- Brazil.

\bibliographystyle{amsalpha}
\bibliography{carolina}

\end{document}